\documentstyle[12pt]{amsart}
\date{  \today}

\newtheorem{th}{Theorem}[section]  
\newtheorem{dfn}[th]{Definition}   \newtheorem{lem}[th]{Lemma}
    \newtheorem{prp}[th]{Proposition}
\newtheorem{conj}[th]{Conjecture}  \newtheorem{ques}[th]{Question}

		\newcommand{\rar}{\rightarrow}
		\newcommand{\bfp}{{\Bbb{P}}}
		\newcommand{\bfq}{{\Bbb{Q}}}

		\newcommand{\co}{{\cal{O}}}
		\newcommand{\das}{\dashrightarrow}

	\newcommand{\Y}{{\cal Y}}

\newcommand{\Kod}{{\operatorname{Kod}}}

\newcommand{\ddim}{\operatorname{Ddim}}

\begin{document}

\noindent %
{\Large LANG MAPS AND HARRIS'S CONJECTURE

\noindent %
{\large\tt PRELIMINARY VERSION} }\\[2mm]
\noindent %
{Dan Abramovich\footnote{Partially supported by NSF grant
DMS-9503276.} \\ 
Department of Mathematics, Boston University\\
111 Cummington, Boston, MA 02215, USA \\
{\tt abrmovic@@math.bu.edu}}\\[2mm]
\today

\noindent %

\large
\addtocounter{section}{-1}
\section{Introduction}
We work over fields of characteristic 0.

Let $X$ be a variety of general type defined over a number field $K$. A well
known conjecture of S. Lang \cite{langbul} states that the set of rational
points $X(K)$ is not Zariski - dense in $X$. As noted in \cite{a}, this implies
that if $X$ is a variety which only {\em dominates} a variety of general type
then $X(K)$ is still not dense in $X$.

J. Harris proposed a way to quantify this situation \cite{harris}: define the
{\em Lang dimension} of a variety to be the maximal dimension of a variety of
general type which it dominates. Harris conjectured in particular that if the
Lang dimension is 0 then for {\em some} number field $L\supset K$ we have that
the set of $L$ rational points $X(L)$ is dense in $X$. The full statement of
Harris's conjecture will be given below (Conjecture \ref{harcon}).

The purpose of this note is to provide a geometric context for Harris's
conjecture, by showing the existence of a universal dominant map to a variety
of general type, which we call {\em the Lang map.}

\section{The Lang map}
\begin{th}
Let $X$ be an irreducible variety ofver a field $k$, $char(k)=0$. There exists
a variety of general type $L(X)$ and a dominant rational map $L:X
\das L(X)$, defined over $k$, satisfying the following universal
property:

 Given a field $K \supset k$ and a dominant rational map $f:X_K\das Z$
defined over $K$, where $Z$ is of general type, then there exists a unique
dominant rational map $L(f): L(X)\das Z$ such that $L(f)\circ L = f$.
\end{th}

\begin{dfn} The universal dominant rational map $L:X\das L(X) $ is
called {\em the Lang map}\footnote{I believe this name is
appropriate since (1) $L$ is closely related to Lang's conjecture, and (2)
the construction resembles in many ways the constructions of Albanese, trace
etc. in Lang's book \cite{langAV}.
I'm also surprised that such a definition has not been made
before. (Is that true?)}. The dimension $\dim L(X)$ is called {\em the Lang
dimension} of $X$.
\end{dfn}

\begin{lem}\label{two-maps} Assume $f_i:X_K\das Z_i$ are dominant rational maps
with
irreducible general fiber, where $Z_i$ are varieties of general type,
$i=1,2$. Then there exists a variety of general type $Z$ over $K$ and dominant
rational maps
$f:X\das Z$ and $g_i: Z \das Z_i$ such that $g_i f = f_i$.
\end{lem}

{\bf Proof}.
 Let $Z = Im(f_1 \times f_2:X \das
Z_1\times Z_2),$ and let $f:X\das Z$ be the induced map. The map $g_i:Z\das
Z_i$ is
dominant and has irreducible general fiber. We claim that $Z$ is a variety of
general type.
By Viehweg's additivity theorem (\cite{viehweg1}, Satz III), it
suffices to show that the generic fiber of $g_1$ is of general type. This
follows since the fibers of $g_1$ sweep $Z_{2}$. (Specifically, let $d$ be the
dimension of the generic fiber of $g_2$. Choose a
general codimension-$d$ plane section  $H\subset Z_1$, then $g_1^{-1}H\rar Z_2$
is generically finite and dominant, therefore $g_1^{-1}H$ is of general type,
therefore the generic fiber of $g_1^{-1}H\rar Z_1$ is of general type.)\qed

\begin{lem} Given a  field extension $K\supset k$, let $l_K$ be the maximal
dimension of a variety of general type  $Y/K$ such that there exists a
dominant rational map $L_K: X_K\das Y$ with irreducible general fiber.
Let $l=\max_{K\supset k}l_K$, and let $K\supset k$ be an extension such that
$l_K=l$.  Then any such map $L_K$ is the Lang map of $X_K$.
\end{lem}

{\bf Proof.} Given an extension $E\supset K$ let $f_2:X_{E} \das Z_2$ be a
dominant rational map.  By the lemma above with $L_K=f_1$  there
exists a variety $Z$ and a dominant rational map with irreducible general fiber
$f:X_E\das Z$ dominating both $Y_E$ and $Z_2$. By maximality $\dim Y= \dim Z_2$
and since the general fibers of $Z\das Y_E$ are irreducible, $Z\das Y_E$ is
birational. The map $g_2 \circ g_1^{-1}:Y_E\das Z_2$ gives the required
dominant rational map. \qed

{\bf Proof of the theorem.}
 Using Stein factorization we may restrict attention
to maps with irreducible general fibers. As above, let $l=\max_{K\supset
k}l_K$,
and let $K\supset k$ be an extension such that $l_K=l$.
We need to show that $L_K$ can be descended to $k$.

 First, we may assume that
$K$ is finitely generated over $k$, since both $Y$ and $L_K$ require only
finitely many coefficients in their defining equations.

Next, we descend $L_K$
to an algebraic extension of $k$. Choose a model $B$ for $K$, and a model
$\Y\rar B$ for $Y$. We have a dominant rational map $X_B\das \Y$ over
$B$. There exists a point $p\in B$  with $[k(p):k]$ finite, such
that $\Y_p$ is a variety of general type of dimension $l$ and such that the
rational map $X_p\das \Y_p$ exists. The lemma above shows that $(\Y_p)_K$ is
birational to $Y$. Alternatively, this step follows
since by theorems of Maehara (see
\cite{moriwaki}) and Kobayashi - Ochiai (see \cite{d-m}) the set of rational
maps to varieties of general type $X\das Z$ is discrete, therefore
each $f:X_K \das Z$ is birationally equivalent to a map defined
over a finite extension of $k$.

We may therefore replace $K$ by an algebraic Galois extension of $k$, which we
still call $K$. Let $Gal(K/k)=\{\sigma_1,\ldots,\sigma_m\}$. For any $1\leq
i\leq m$ we have a rational map
$(f_1\times\sigma_i\circ f_1): X\das Y\times Y^{\sigma_i}$. Applying lemma
\ref{two-maps} we obtain a birational map
$Y\das
Y^{\sigma}.$ There are open sets $U_i\subset Y^{\sigma_i}$ over which these
maps are regular isomorphisms, giving rise to descent data for $U_1$ to
$k$. \qed

Is there a way to describe the fibers of the Lang map $X\das L(X)$? A first
approximation is provided by the following:

\begin{prp} The generic fiber of the Lang map has Lang dimension 0.
\end{prp}

{\bf Proof. }  Let $\eta\in L(X)$ be the generic point and let $X_\eta \das
L(X_\eta)$ be the lang map of the generic fiber. Let $M\rar L(X)$ be a model of
$L(X_\eta)$. By definition, the generic fiber $M_\eta=L(X_\eta)$ of $M$ is of
general type, therefore by Viehweg's additivity theorem $M$ is of general type,
and by definition $M$ is birational to $L(X)$.\qed

\begin{ques}\label{open}
Is there an open set in $X$ where the Lang map is defined and the fibers have
Lang dimension 0?
\end{ques}

 We will see that the answer is yes, if one assumes the following inspiring
conjecture of higher dimensional classification theory:

\begin{conj}[see Conjecture 1.24 of \cite{flab}]\label{uni/kod}
\begin{enumerate} \item Let $X$ be a
variety in characteristic 0. Then either $X$ is uniruled, or $\Kod(X) \geq 0$.
\item If $\Kod(X) \geq 0$ then there is an open set in $X$ where the fibers of
the Iitaka fibration have Kodaira dimension 0.
\end{enumerate}
\end{conj}

This conjecture allows us to ``construct'' the Lang map ``from above'':

\begin{prp}\label{constr} Assume that conjecture \ref{uni/kod} holds true. Then
there is a
finite sequence of dominant rational maps
$$X\das X_1\das\cdots \das X_n = L(X)$$
where each map $X_i\das X_{i+1}$ is either an MRC fibration (see \cite{kmm},
2.7) or an Iitaka fibration. In particular, the answer to question \ref{open}
is ``yes''.
\end{prp}

{\bf The proof} is obvious. We remark that since \ref{uni/kod} is known when
the fibers have dimension $\leq 2$. In particular, \ref{constr} is known
unconditionally  when $\dim X\leq 3$.

\section{Harris's conjecture}
As mentioned above, we define {\em the Lang dimension} of a variety $X$ to be
$\dim L(X)$, and Lang's conjecture implies that if $K$ is a number field, and
if $X/K$ has positive Lang dimension, then $X(K)$ is not Zariski - dense in
$X$. In \cite{harris}, J. Harris proposed a complementary statement:

\begin{conj}[Harris's conjecture, weak form]\label{harwk} Let $X$ be a variety
of Lang dimension 0
defined over a number field $K$. Then for some finite extension $E\supset K$
the set of $E$-rational points $X(E)$ is Zariski dense in $X$.
\end{conj}
 It is illuminating to consider the motivating case of an elliptic
surface of positive rank.

Let $\pi_0:X_0\rar \bfp^1$ be a pencil of cubics through 9 rational points in
$\bfp^2$. By choosing the base points in general position we can guarantee that
the pencil has 12 irreducible singular fibers which are nodal rational
curves. The Mordell - Weil group of $\pi_0$ has rank 8. The relative dualizing
sheaf $\omega_{\pi_0}= \co_{X_0}(F_0)$ where $F_0$ is a fiber. Let
$f:\bfp^1\rar\bfp^1$ be a map of degree at least 3. Let $\pi:X\rar\bfp^1$ be
the pull-back of $X_0$ along $f$. Then $\omega_\pi = \co_X(3F)$, therefore
$\omega_X = \co_X(F)$ and $X$ has Kodaira dimension 1. The Iitaka fibration is
simply $\pi$. The elliptic surface $X$ still has a Mordell - Weil group of rank
8 of sections. By applying these sections to rational points on $\bfp^1$ we see
that the set of rational points $X(\bfq)$ is dense in $X$.

It is not hard to modify this example to obtain a varying family of elliptic
surfaces which has a dense collection of sections. Let $B$ be a curve and let
$g: B\times \bfp^1 \rar \bfp^1$ be a family of rational functions on $\bfp^1$
which varies in moduli (such families exists as soon as the degree is at least
3). Let $Y$ be the pullback of $X$ to $B\times \bfp^1$. Then $p: Y\rar B$ is a
family of elliptic surfaces,  of variation $Var(p)=1$, and relative kodaira
dimension 1. By composing sections of $E$ with $g$ and arbitrary rational maps
$B\rar \bfp^1$, we see that $p$ has a dense collection of sections.

Harris's weak conjecture for elliptic surfaces is attributed to
Manin.  In case of surfaces over $\bfp^1_\bfq$, it has been related to
the conjecture of Birch and Swinnerton-Dyer: let $\pi:X\rar \bfp^1$ be an
elliptic surface defined over $\bfq$. In \cite{manduchi}, E. Manduchi
shows that under certain assumptions on the behavior of the $j$ function, the
set of points in $\bfp^1(\bfq)$ where the fiber has root number $-1$ is dense
(in the classical topology). According to the conjecture of Birch and
Swinnerton-Dyer, the root number gives the parity of the Mordell-Weil rank. It
is likely that some of Manduchi's conditions (at least condition (1) in Theorem
1  of \cite{manduchi})  can be relaxed once one passes to a number field.

What can be said in case $0<\dim L(X) <\dim X$? In \cite{harris2}, Harris
proposed the
following definition:

\begin{dfn} The {\em diophantine dimension,} $\ddim(X)$ is defined as follows:
$$\begin{array}{cccc} \ddim(X) :=& \min & \max & \dim(\overline{U(E)})\\[-2mm]
 & _{\emptyset \neq U\subset X\mbox{ open   }} & _{[E:K]<\infty} &
\end{array}$$
\end{dfn}

Harris proceeded to propose the following:

\begin{conj}[Harris's conjecture]\label{harcon} For any variety $X$ over
a number field, $$ \ddim(X) + \dim L(X) = \dim X.$$
\end{conj}

I do not know whether or not Harris himself believes this conjecture. This does
not really matter. What is appealing in this conjecture, apart from it's
``tightness'', is that any evidence, either for or against it, is likely to
be of much interest.

For lack of any better results, we just note that proposition \ref{constr}
directly implies the following:

\begin{prp} Assuming \ref{uni/kod}, Lang's conjecture together with the weak
form of Harris's conjecture \ref{harwk} implies Harris's conjecture
\ref{harcon}.
\end{prp}

{\sc Acknowledgements } \small I would like to thank D. Bertrand, J. Harris,
J. Koll\'ar, K. Matsuki, D. Rohrlich and J. F. Voloch for discussions related
to this note.

\end{document}